\documentclass[12pt]{article}

\usepackage[margin=2.54cm]{geometry}
\usepackage{setspace}
\usepackage{amscd,amssymb, amsmath, setspace, graphicx}
\usepackage{url, verbatim}
\usepackage{hyperref}
\usepackage{color}

\usepackage{tikz}
\usetikzlibrary{shapes,arrows}
\tikzstyle{decision} = [diamond, draw, fill=blue!20,
    text width=4.5em, text badly centered, node distance=2.5cm, inner sep=0pt]
\tikzstyle{longblock} = [rectangle, draw, fill=blue!20,
    text width=15em, text centered, rounded corners, minimum height=2.5em]
\tikzstyle{block} = [rectangle, draw, fill=blue!20,
    text width=4.5em, text centered, rounded corners, minimum height=2.5em]
\tikzstyle{line} = [draw, very thick, color=black!50, -latex']
\tikzstyle{cloud} = [draw, ellipse,fill=red!20, node distance=2.5cm,
    minimum height=2em]

\newcommand{\macaulay}{\mbox{\sc Macaulay2}}
\newcommand{\bertini}{\mbox{\sc Bertini}}
\newcommand{\bertiniM}{\mbox{\it Bertini}}
\newcommand{\NAGtypes}{\mbox{\it NAGtypes}}
\newcommand{\NAG}{\mbox{\it NumericalAlgebraicGeometry}}
\newcommand{\PHCpack}{\mbox{\sc PHCpack}}

\title{Bertini for Macaulay2}

\author{Daniel J. Bates\thanks{Department of Mathematics, Colorado State University, Fort Collins, CO 80523. {\tt bates@math.colostate.edu}.  This work was partially supported by NSF grant DMS--1115668.}, 
Elizabeth Gross\thanks{Department of Mathematics, North Carolina State University, Raleigh, NC 27695; {\tt eagross@ncsu.edu}. This work was partially supported by NSF award DMS--1304167.}, 
Anton Leykin\thanks{School of Mathematics, Georgia Tech, Atlanta, GA; {\tt leykin@gatech.math.edu}. This work was partially supported by NSF grant DMS--1151297.}, 
Jose Israel Rodriguez\thanks{Department of Mathematics, University of California at Berkeley, 
 Berkeley, CA 94720 
 {\tt jo.ro@berkeley.edu}.
This work was partially supported by NSF grant DMS--0943745.}}

\begin{document}

\date{\today}

\maketitle

\begin{abstract} 
Numerical algebraic geometry is the field of computational mathematics 
concerning the numerical solution of polynomial systems of equations.  
\bertini, a popular software package for computational applications of this field,
includes implementations of a variety of algorithms based on polynomial homotopy continuation.  The \macaulay\ package \bertiniM\ provides an interface to \bertini, making it possible to access the core run modes of \bertini\ in \macaulay.  With these run modes, users can find approximate solutions to zero-dimensional systems and positive-dimensional systems, test numerically whether a point lies on a variety, sample numerically from a variety, and perform parameter homotopy runs.
\end{abstract}

\section{Numerical algebraic geometry }

{\em Numerical algebraic geometry (numerical AG)} refers to a set of methods for finding and manipulating the 
solution sets of systems of polynomial equations.  Said differently, given $f:\mathbb C^N\to\mathbb C^n$, 
numerical algebraic geometry provides facilities for computing numerical approximations to isolated 
solutions of $V(f)=\left\{z\in\mathbb C^N | f(z)=0\right\}$, as well as numerical approximations to 
generic points on positive-dimensional components.  The book~\cite{SWbook} provides a good 
introduction to the field, while the newer book~\cite{BertiniBook} provides a simpler introduction 
as well as a complete manual for the software package \bertini~\cite{Bertini}.

\bertini\ is a free, open source software package for computations in numerical algebraic geometry.
The purpose of this article is to present a \macaulay\ \cite{M2} package \bertiniM\ that provides an interface to \bertini.  
This package uses basic datatypes and service routines for computations in numerical AG provided by the package \NAGtypes. 
It also fits the framework of the \NAG\ package~\cite{NAG4M2}, a native \macaulay\ implementation of a collection of numerical AG algorithms: most of the core functions of \NAG\ have an option of using \bertini\ instead of the native solver.

In the remainder of this section, we very briefly describe a few fundamental concepts 
of the field.  In the subsequent sections, we describe the various run modes of \bertini\ that have 
been implemented in this interface.  We conclude with Section 5, which describes how to use \bertiniM\ within \NAG.

\subsection{Finding isolated solutions}

The core computational engine within \bertini\ is {\em homotopy continuation}.  This is a 
three-stage process for finding a superset of all isolated solutions in $V(f)$.  Given a polynomial 
system $f(z)$, the three steps are as follows:
\begin{enumerate}
\item Choose an easily-solved polynomial system $g(z)$ that reflects the structure of $f(z)$, 
and solve it.  Call this set of solutions $S$.
\item Form the homotopy 
$$H(z,t)=(1-t)f(z) + \gamma tg(z),$$
with $\gamma\in\mathbb C$ a random complex number.  Notice that $H(z,1)=\gamma g(z)$, the 
solutions of which are known, and $H(z,0)=f(z)$, for which we seek the solutions.
\item There is a real curve extending from each solution $z\in S$.  Use predictor-corrector methods, 
adaptive precision, and endgames to track along all of these paths as $t$ goes from 1 to 0.
\end{enumerate}
Assuming $g(z)$ is constructed in one of several canonical ways~\cite{SWbook}, there is a 
probability one guarantee that this procedure will result in a superset of all isolated solutions 
of $f(z)=0$.

There are many variations on this general technique, and there are many minor issues to 
consider when implementing this method.  However, due to space limitations, we leave the reader 
to explore the references for more information on this powerful method.

\subsection{Finding irreducible components}

Given an irreducible algebraic set $X$ of dimension $k$, it is well known that $X$ will intersect 
almost any linear space of codimension $k$ in a finite set of points.  In fact, there is a Zariski open 
subset of the set of all linear spaces of codimension $k$ for which intersection with $X$ yields 
some fixed number of points, called the {\em degree} of $X$, $\deg X$.

This fundamental fact underlies the computation of positive-dimensional irreducible components 
in numerical algebraic geometry.  Suppose algebraic set $Z$ decomposes into irreducible 
components $Z_{i,j}$, 
$$Z=\cup_{i=0}^{\dim Z}\cup_{j\in \Lambda_i} Z_{i,j},$$
where $i$ is the dimension of $Z_{i,j}$ and $j$ is just the index of component $Z_{i,j}$ in dimension 
$i$, stored in finite indexing set $\Lambda_i$.

In numerical algebraic geometry, the representation $W$ of an algebraic set $Z$ consists of a representation 
$W_{i,j}$ for each irreducible component $Z_{i,j}$ of $Z$.  In particular, {\em witness set} $W_{i,j}$ is 
a triple $(f, L_{i}, \widehat{W}_{i,j})$, consisting of polynomial system $f$, linear functions $L_{i}$ 
corresponding to a linear space of  codimension $i$, and {\em witness point set} 
$\widehat{W}_{i,j} = Z_{i,j}\cap V(L_{i})$.  

There are a variety of ways to compute $W$, many of which are described in detail in~\cite{BertiniBook}.  Most of these methods can be accessed through the package \bertiniM\ by using optional inputs to specify the desired algorithm.

\section{Solving zero-dimensional systems}

In the following sections we outline and give examples of the different \bertini\ run modes implemented in the interface package \bertiniM.

\subsection{Finding solutions of zero-dimenstional systems}

The method {\tt bertiniZeroDimSolve} calls \bertini\ \ to solve a
polynomial system and returns solutions as a list of {\tt Points} using the data types from  \NAGtypes. Diagnostic information, such as the residuals and the condition number, are stored with the coordinates of the solution and can be viewed using {\tt peek}.

\begin{verbatim}
i1 : R=CC[x,y];
i2 : f = {x^2+y^2-1,(x-1)^2+y^2-1};
i3 : solutions=bertiniZeroDimSolve(f)
o3 = {{.5, .866025}, {.5, -.866025}}
i4 : peek solutions_0
o4 = Point{ConditionNumber => 88.2015     }
            Coordinates => {.5, .866025}
            CycleNumber => 1
            FunctionResidual => 3.66205e-15
            LastT => .000390625
            MaximumPrecision => 52
            NewtonResidual => 4.27908e-15
            SolutionNumber => 3
\end{verbatim}

Users can specify to use regeneration, an equation-by-equation solving method, by setting the option {\tt USEREGENERATION} to 1.

\begin{verbatim}
i5 : solutions=bertiniZeroDimSolve(f, USEREGENERATION=>1);
\end{verbatim}

In common applications, one would like to classify solutions, e.g., separate real solutions from non-real solutions, and, thus, recomputing solutions to a higher accuracy becomes important.  The method {\tt bertiniRefineSols} calls the sharpening module of \bertini\  and sharpens a list of solutions to a desired number of digits using Newton's method.

\begin{verbatim}
i6 : refinedSols=bertiniRefineSols(f, solutions, 20);
     (coordinates refinedSols_0)_1
o6 = .86602540378443859659+3.5796948761134507351e-83*ii
\end{verbatim}     

\subsection{Parameter homotopies}

Many fields, such as statistics, physics, biochemistry, and engineering, have applications that require solving a large number of systems from a parameterized family of polynomial systems.  In such situations, computational time can be decreased by using parameter homotopies.  For an example illustrating how parameter homotopies can be used in statistics see \cite{HRS12}. 

The method {\tt bertiniParameterHomotopy} calls \bertini\ to run both stages of a parameter homotopy.  First, \bertini\  assigns a random complex number to each specified parameter and solves the resulting system, then, after this initial phase, \bertini\ computes solutions for every given choice of parameters using a number of paths equal to the exact root count.

\begin{verbatim}
i7 : R=CC[a,b,c][x,y];
i8 : f={a*x^2+b*y^2-c, y};
i9 : bertiniParameterHomotopy(f,{a,b,c},{{1,1,1},{2,3,4}})
o9 = {{{-1, 0}, {1, 0}}, {{-1.41421, 0}, {1.41421, 0}}}
\end{verbatim}

\section{Solving positive-dimensional systems}

Given a positive-dimensional system $f$, the method {\tt bertiniPosDimSolve} calls \bertini\ to compute a numerical irreducible decomposition.  This decomposition is assigned the type {\tt NumericalVariety} in \macaulay.  In the default settings, \bertini\ uses a classical cascade homotopy to find witness supersets in each dimension, removes extra points using a membership test or local dimension test, deflates singular witness points, then factors using a combination of monodromy and a linear trace test.

\begin{verbatim}
i10 :    R = CC[x,y,z];
i11 :    f = {(y^2+x^2+z^2-1)*x, (y^2+x^2+z^2-1)*y};
i12 :    NV = bertiniPosDimSolve f

o12 = A variety of dimension 2 with components in
      dim 1:  [dim=1,deg=1]
      dim 2:  [dim=2,deg=2]. 

o12 : NumericalVariety
\end{verbatim}

Once the solution set to a system, i.e., the variety $V$, is computed and stored as  a {\tt NumericalVariety},  {\tt bertiniComponentMemberTest} can be used to test numerically whether a set of points $p$ lie on the variety $V$. For every point in $p$, {\tt bertiniComponentMemberTest} returns the components to which that point belongs. As for sampling, {\tt bertiniSample}  will sample from a witness set $W$. These methods call the membership testing and sampling options in \bertini\ respectively.

\begin{verbatim}
i13 : p={{0,0,0}};
i14 : bertiniComponentMemberTest (NV, p)
o14 = {{[dim=1,deg=1]}}

i15 : component=NV#1_0
i16 : bertiniSample(component,1)
o16 = {{0, -8.49385e-20+7.48874e-20*ii, -.148227-.269849*ii}}
\end{verbatim}

\section{Solving homogeneous systems}
The package \bertiniM\ includes functionality to solve a homogenous system  that defines a projective variety. In \bertini, the numerical computations are performed on a generic affine chart to  compute representatives of  projective points.  To solve homogeneous  equations,  set the option $\tt ISPROJECTIVE$ to  $1$. 
If the user  inputs a square system of $n$ homogeneous equations in $n+1$ unknowns, then  the method {\tt bertiniZeroDimSolve} 
outputs a list of points in projective space. 

\begin{verbatim}
i17 : R = CC[x,y,z];
i18 : f = {y^2-4*z^2,16*x^2-y^2};
i19 : bertiniZeroDimSolve(f,ISPROJECTIVE=>1);
o19 = {{.251411+.456072*ii, 1.00564+1.82429*ii, .502821+.912143*ii}, 
 {.106019+.160896*ii, .424078+.643585*ii, -.212039-.321792*ii},
 {-.15916-.12286*ii, .636639+.49144*ii, -.318319-.24572*ii},
 {-.48005-.092532*ii, 1.9202+.370128*ii, .960101+.185064*ii}}
\end{verbatim} 

If $f$ is a positive-dimensional homogeneous system of equations, then the method \break
 {\tt bertiniPosDimSolve} calls \bertini\ to compute a numerical irreducible decomposition of the projective variety defined by $f$.

\begin{verbatim}
i20 : R = CC[x,y,z];
i21 : f = {(x^2+y^2-z^2)*(z-x),(x^2+y^2-z^2)*(z+y)};
i22 : NV = bertiniPosDimSolve(f,ISPROJECTIVE=>1)
o22 : = A projective variety with components in projective dimension: 
      dim 0:  [dim=0,deg=1]
      dim 1:  [dim=1,deg=2]
\end{verbatim}

\section{Using \bertiniM\ from \NAG}
The \bertiniM\ package depends on the \NAGtypes\ package, a collection of basic datatypes and service routines common to all \macaulay\ packages for numerical AG: e.g., an interface package~\cite{PHCpackM2} to another polynomial homotopy continuation solver, \PHCpack~\cite{PHC}, also has this dependence. The dependencies between the metioned packages are depicted on the diagram below: dashed arrows stand for a dependency that is optional and is engaged if an executable for the corresponding software is installed.

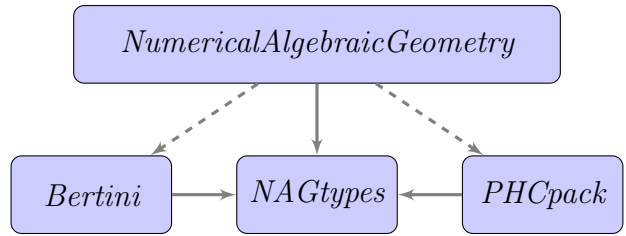
\begin{center}
\begin{tikzpicture}[scale=2, node distance = 2cm, auto]
    \node [longblock] (NAG) {{\em NumericalAlgebraicGeometry}};
    \node [block, below of=NAG] (NAGtypes) {{\em NAGtypes}};
    \node [block, left of=NAGtypes, node distance=3cm] (Bertini) {{\em Bertini}};
    \node [block, right of=NAGtypes, node distance=3cm] (PHCpack) {{\em PHCpack}};
    \path [line] (NAG) -- (NAGtypes);
    \path [line] (Bertini) -- (NAGtypes);
    \path [line] (PHCpack) -- (NAGtypes);
    \path [line,dashed] (NAG) -- (Bertini);
    \path [line,dashed] (NAG) -- (PHCpack);
\end{tikzpicture}
\end{center}

While independent from the \NAG\ package, our interface provides a valuable option for this package: the user can set \bertini\ as a default solver for homotopy continuation tasks.
\begin{verbatim}
i23 : needsPackage "NumericalAlgebraicGeometry";
i24 : setDefault(Software=>BERTINI)
\end{verbatim}
An alternative way is to specify the {\tt Software} option in a particular command:
\begin{verbatim} 
i25 : CC[x,y]; system = {x^2+y^2-1,x-y};
i26 : sols = solveSystem(system, Software=>M2engine)
o26 = {{-.707107, -.707107}, {.707107, .707107}}
i27 : refsols = refine(system, sols, Bits=>99, Software=>BERTINI);
i28 : first coordinates first refsols
o28 = -.707106781186547524400844362105-2.13764004134262114934289955768e-140*ii
o28 : CC (of precision 100)
\end{verbatim}
The unified framework for various implementations of numerical AG algorithms should be particularly convenient to a \macaulay\ user doing numerical computations with tools from many packages.


\begin{thebibliography}{99}

\bibitem{Bertini}
D.J. Bates, J.D. Hauenstein, A.J. Sommese, and C.W. Wampler.
Bertini: Software for numerical algebraic geometry.
Available at {\tt http://www.nd.edu/$\sim$sommese/bertini}.

\bibitem{BertiniBook}
D.J. Bates, J.D. Hauenstein, A.J. Sommese, and C.W. Wampler.  
{\em Numerically solving polynomial systems with Bertini.}
To be published by SIAM, 2013.

\bibitem{M2}
D.R. Grayson and M.E. Stillman. Macaulay2, a software system for research in algebraic geometry, available at {\tt www.math.uiuc.edu/Macaulay2/.}

\bibitem{PHCpackM2}
E. Gross, S. Petrovi\'{c}, and J. Verschelde.
{\em Interfacing with PHCpack.}
J. Softw. Algebra Geom. 5, 20--25, 2013.

\bibitem{HRS12} J. Hauenstein, J. Rodriguez, and B. Sturmfels. \emph{Maximum likelihood for matrices with rank constraints}, preprint, arXiv:1210.0198.

\bibitem{NAG4M2}
A. Leykin. 
{\em Numerical algebraic geometry.}
J. Softw. Algebra Geom. 3, 5--10, 2011.

\bibitem{SWbook}
A.J. Sommese and C.W. Wampler.
{\em The numerical solution to systems of polynomials arising in engineering and science.}
World Scientific, Singapore, 2005.

\bibitem{PHC}
J. Verschelde.
{\em Algorithm 795: PHCpack: A general-purpose solver for polynomial systems by homotopy continuation} 
ACM Trans. Math. Softw. 25(2):251--276, 1999.  Software available at
{\tt http://www.math.uic.edu/$\sim$jan}.


\end{thebibliography}
\end{document}